\documentclass[12pt]{article}

\usepackage{amssymb,amsthm,amsmath,eucal,latexsym}
\usepackage{color,caption,comment}
\usepackage{multirow}
\usepackage{listings}
\usepackage{setspace}

\usepackage{changepage}
\usepackage[linesnumbered,lined,boxed,commentsnumbered]{algorithm2e}
\usepackage{verbatim}
\usepackage{tikz}
\usepackage{xcolor}
\usepackage{algpseudocode}

\newtheorem{theorem}{Theorem}[section]
\newtheorem{corollary}[theorem]{Corollary}
\newtheorem{lemma}[theorem]{Lemma}
\newtheorem{conjecture}[theorem]{Conjecture}
\newtheorem{definition}[theorem]{Definition}
\newtheorem{remark}[theorem]{Remark}

\usepackage[utf8]{inputenc}

\usepackage{enumitem, graphicx, hyperref}
\graphicspath{ {./XXfiguresXX/} }
\allowdisplaybreaks

\title{Proportion-Based Hypergraph Burning}
\author{Andrea C. Burgess\thanks{Department of Mathematics and Statistics, University of New Brunswick, Saint John, NB, E2L~4L5, Canada. \texttt{andrea.burgess@unb.ca}}, John A. Hawkin\thanks{Verafin, a NASDAQ Company, St. John's, NL, A1A~0L9, Canada. \texttt{John.Hawkin@verafin.com}}, Alexander J. M. Howse\thanks{Verafin, a NASDAQ Company, St. John's, NL, A1A~0L9, Canada. \texttt{alexander.howse@verafin.com}},\\ Caleb W. Jones\thanks{Department of Mathematics, Toronto Metropolitan University, Toronto, ON, M5B~2K3, Canada. \texttt{caleb.w.jones@torontomu.ca}}, David A. Pike\thanks{Department of Mathematics and Statistics, Memorial University of Newfoundland, St.~John's, NL, A1C~5S7, Canada. \texttt{dapike@mun.ca}}}

\definecolor{darkgreen}{RGB}{34, 175, 65}

\begin{document}

\date{\today}
\maketitle

\begin{abstract}
Graph burning is a discrete process that models the spread of influence through a network using a fire as a proxy for the type of influence being spread. This process was recently extended to hypergraphs. We introduce a variant of hypergraph burning that uses an alternative propagation rule for how the fire spreads -- if some fixed proportion of vertices are on fire in a hyperedge, then in the next round the entire hyperedge catches fire. This new variant has more potential for applications than the original model, and it is similarly viable for obtaining deep theoretical results. We obtain bounds which apply to general hypergraphs, and introduce the concept of the \emph{burning distribution}, which describes how the model  changes as the proportion ranges over $(0,1)$. We also obtain computational results which suggest there is a strong correlation between the automorphism group order and the lazy burning number of a balanced incomplete block design.
\end{abstract}

\section{Introduction} 
\label{introduction}

Graph burning is a game played on finite, simple, undirected graphs over a discrete sequence of rounds. There is a single player, called the \textit{arsonist}, who attempts to set fire to every vertex of the graph in as short a time as possible. 
In each round the arsonist sets a vertex on fire, and existing fires spread or {\em propagate} from burning vertices to adjacent unburned vertices. 
Once a vertex is on fire, it will remain that way until the end of the game. The game ends when each vertex in the graph is on fire.

Graph burning was introduced by Bonato, Janssen and Roshanbin~\cite{first_paper}, although a problem equivalent to the special case of graph burning on an $n$-cube had previously been posed by Brandenburg and Scott~\cite{VERY_first_ref}. For a survey of known results on graph burning, we refer the reader to~\cite{Bonato_summary}, as well as~\cite{new_paper, asymptotically_paper} for recent results.


We now describe the rules for graph burning in detail. Rounds are indexed by positive integers. 
Denote the set of vertices that are on fire at the end of round $r$ by $F_r$, and set $F_0=\emptyset$. During each round $r\geq 1$, the following two things happen simultaneously. 

\begin{itemize}
\item For each $u \in F_{r-1}$ and $v \in V(G) \setminus F_{r-1}$, if $uv \in E(G)$, then $v$ catches fire {\em by propagation}. 
\item The arsonist chooses a vertex $u_r \in V(G)\setminus F_{r-1}$ and sets it on fire. The vertex $u_r$ is called a \textit{source}.
\end{itemize}

In round 1, no vertices catch fire due to propagation, and the arsonist chooses the first source. In round $r>1$, the arsonist is permitted to choose a source which is adjacent to a vertex of $F_{r-1}$, which is already on fire.  We call such sources {\em redundant}.  A redundant source would catch fire in round $r$ without the arsonist's intervention, so it is never advantageous for the arsonist to choose a redundant source if there is another choice.  However, choosing a redundant source may be unavoidable in the final round.

A source is called \emph{valid} if 
it is the $i^{th}$ source chosen, and it is not on fire at the end of round $i-1$. 
A sequence of valid sources, indexed by the round at which they are set on fire, is also called \textit{valid}. Note that redundant sources are also valid. 
A valid sequence of sources $(u_1,u_2,\ldots,u_k)$ that leaves the graph completely burned when the arsonist burns $u_i$ in round $i$ is called a \textit{burning sequence}. 
Given a graph $G$, the {\em burning number} of $G$, denoted $b(G)$, is the minimum length of a burning sequence in $G$.  We call a burning sequence of length $b(G)$ {\em optimal}. Note that $b(G)$ may equivalently be defined as the earliest round at which fire could possibly be spread to all vertices of $G$.

We now give a brief overview of some concepts from hypergraph theory.
A \emph{hypergraph} $H$ is an ordered pair $(V,E)$ where $V$ and $E$ are disjoint finite sets, $V\neq\emptyset$, and each element of $E$ is a subset of $V$. The elements of $V=V(H)$ are called \emph{vertices}, and the elements of $E=E(H)$ are called \emph{edges} or \emph{hyperedges}. Informally, a hypergraph is 
a generalization of a graph, such that edges can now contain any number of vertices, not just two. 
If two vertices belong to a common edge, they are called \emph{adjacent}. The definitions for paths and connectedness/disconnectedness are analogous to those from graph theory. If a vertex does not belong to any edge, it is called \emph{isolated}. A hypergraph is $k$\emph{-uniform} if every edge has cardinality $k$.  Note that a simple graph is a $2$-uniform hypergraph.
Two edges $e_1$ and $e_2$ are \emph{parallel} if they contain exactly the same vertices. The number of edges parallel to some edge $e$, including $e$ itself, is the \emph{multiplicity} of $e$. A hypergraph is called \emph{simple} if no edge has multiplicity greater than 1, and no edge contains at most 1 vertex.
For further information on hypergraphs, see~\cite{bahm_sanj}, \cite{ hypergraph_text}, and \cite{conn_in_hyp}.

\emph{Hypergraph burning} was introduced in \cite{our_paper!, mythesis}  as a natural extension of graph burning. It was later studied on Steiner triple systems \cite{our_STS_paper}  and Latin squares \cite{our_Latin_paper}. The rules for hypergraph burning are identical to those from graph burning, save for the rule for how the fire propagates within a hyperedge:  fire spreads to a vertex $v$ in round $r$ if and only if there is an edge $\{v,u_1,\ldots,u_k\}$ such that each of $u_1,u_2,\ldots,u_k$ was on fire at the end of round $r-1$. Clearly this reduces to the original game of graph burning when the hypergraph has only edges of size 2. The definitions of a \textit{source}, \textit{redundant source}, \textit{(optimal) burning sequence}, and \textit{burning number} are all analogous to those in graph burning. 


There is also an alternative rule set for hypergraph burning, called \emph{lazy hypergraph burning} \cite{our_paper!, mythesis}.  
In this case, a lazy arsonist sets fire to a set of vertices simultaneously. The arsonist's aim is that the remaining vertices of the hypergraph are eventually all burned through subsequent propagation.  A set of initial burning vertices which propagate to the entire hypergraph is called a {\em lazy burning set}, and the {\em lazy burning number} of a hypergraph $H$, denoted $b_L(H)$, is the minimum size of a lazy burning set.  We call a lazy burning set of size $b_L(H)$ {\em optimal}.
For this model, we no longer care about the number of time steps required for the hypergraph to become fully burned, just the size of a smallest lazy burning set. Note that lazy hypergraph burning with the original propagation rule is equivalent to a process known as $\mathcal{H}$\emph{-bootstrap percolation}, which was introduced in \cite{bootstrap_paper_1}. 

Broadly speaking, graph burning is a model for the spread of influence throughout a network. One analogy that is often used is social contagion, in which the fire represents a piece of misinformation, and the graph represents a social network (i.e.~people are represented by vertices and if two people are friends then their corresponding vertices are adjacent). Hence, the burning number is a measure of how quickly everyone in the social network could become misinformed. Hypergraph burning has the potential to improve this model, as misinformation is not only spread via one-on-one friendships in the modern world. In particular, hyperedges of size greater than two could represent social media friend groups. Unfortunately, the original propagation rule for hypergraph burning is not very conducive to the way information spreads online in the real world. We therefore wish to formulate a new propagation rule which is better suited to real-world applications. Hence, the topic of this paper is a variant of (lazy) hypergraph burning that uses a new propagation rule: \\

\begin{centering}
\fbox{\parbox{0.983\textwidth}{For some fixed \emph{proportion} $p\in (0,1)$, if $\left\lceil p |e| \right\rceil$ or more vertices are on fire in an edge $e$, then in the next round all unburned vertices in $e$ catch fire. }}
\end{centering}\\

This rule is much more conducive to the social contagion analogy -- if a large enough proportion of the members of a social media friend group believe a piece of misinformation, it makes sense that they would then spread this misinformation to other members of the group. We will call this variant \emph{proportion-based (lazy) hypergraph burning}. The rules of the new variant are identical to those of (lazy) hypergaph burning, save for the new propagation rule. The definitions of a \textit{(redundant/valid) source}, \textit{(valid/optimal) burning sequence}, and \textit{(optimal) lazy burning set} are all analogous to those in hypergraph burning. 

\begin{definition}
Given a hypergraph $H$ and a proportion $p\in(0,1)$, the \emph{burning number of} $H$ \emph{with respect to} $p$, denoted $b_p(H)$, is the length of an optimal burning sequence for $H$ when using the proportion-based propagation rule with proportion $p$.
\end{definition}

\begin{definition}
Given a hypergraph $H$ and a proportion $p\in(0,1)$, the \emph{lazy burning number of} $H$ \emph{with respect to} $p$, denoted $b_{L,p}(H)$, is the size of an optimal lazy burning set for $H$ when using the proportion-based propagation rule with proportion $p$.
\end{definition}

For the remainder of this paper we will assume that we are using this new propagation rule, unless it is otherwise stated. When we refer to the burning number, we mean $b_p(H)$, and when we refer to the lazy burning number, we mean $b_{L,p}(H)$, unless otherwise stated.

\section{General Results and Bounds}
\label{general_results}

Notice that under the proportion-based propagation rule, if $p$ is high enough, there may be edges in the hypergraph that cannot possibly cause fire to spread. These edges are the subject of our first few results.

\begin{remark}
Let $p\in (0,1)$, $H$ be a hypergraph, and $e\in E(H)$. If $\left\lceil p|e|\right\rceil=|e|$ then $e$ cannot cause fire to propagate among its vertices.
\end{remark}

This observation motivates the following definition.

\begin{definition}
\label{non_flam}
Let  $p\in(0,1)$, $H$ be a hypergraph and $e\in E(H)$. If $\left\lceil p|e|\right\rceil=|e|$ then $e$ is \emph{non-flammable with respect to} $p$. Otherwise, $e$ is \emph{flammable with respect to} $p$. 
\end{definition}

We will often omit the words ``with respect to $p$'' when $p$ is obvious. It is straightforward to determine whether or not an edge is non-flammable.

\begin{lemma}
\label{non_flam_iff}
An edge $e$ is non-flammable if and only if $|e|<\frac{1}{1-p}$.
\end{lemma}

\begin{proof}
Observe that $\left\lceil p|e|\right\rceil=|e|$ if and only if $|e|-p|e|<1$. Keeping in mind that $1-p$ is positive, we can simplify to obtain $|e|<\frac{1}{1-p}$.
\end{proof}

\begin{lemma}
\label{delete_non_flam_edges}
Let  $p\in(0,1)$, $H$ be a hypergraph, and $e_1,e_2,\ldots,e_k$ be all of the non-flammable edges in $H$ (for this value of $p$). Let $H^\prime$ be the hypergraph with vertex set $V(H)$ and edge set $E(H)\setminus\{e_1,e_2,\ldots,e_k\}$. Then $b_p(H)=b_p(H^\prime)$ and $b_{L,p}(H)=b_{L,p}(H^\prime)$.
\end{lemma}

\begin{proof}
Non-flammable edges have no effect on the burning game. Thus, in the context of the burning game, $H$ and $H^\prime$ behave identically.
\end{proof}




In a finite hypergraph we can always choose $p$ close enough to $1$ such that every edge is non-flammable, in which case the following remark applies.

\begin{remark}
\label{all_non_flam_cor}
Let $p\in(0,1)$. If every edge in $H$ is non-flammable (for this value of $p$), then $b_{L,p}(H)=b_p(H)=|V(H)|$.
\end{remark}

Similarly, we can always choose $p$ in such a way as to ensure each non-singleton edge is flammable.

\begin{lemma}
Let $p\in (0,\frac{1}{2}]$ and $H$ be a hypergraph. Then the only non-flammable edges in $H$ are singleton edges.
\end{lemma}

\begin{proof}
Since $0<p\leq\frac{1}{2}$ we can calculate $1<\frac{1}{1-p}\leq 2$. Now, for this value of $p$, an edge $e$ is non-flammable if and only if $|e|<\frac{1}{1-p}$, that is, if $|e|<2$.
\end{proof}

If the proportion $p$ is very low, then there may exist edges $e$ in the hypergraph such that, if a single vertex in $e$ is on fire, then all of $e$ catches fire.

\begin{remark}
Let $p\in (0,1)$, $H$ be a hypergraph, and $e\in E(H)$. If $\left\lceil p|e|\right\rceil=1$ then a single burning vertex in $e$ will cause fire to propagate among all the vertices in $e$.
\end{remark}

This observation motivates the following definition.

\begin{definition}
\label{high_flam}
Let  $p\in(0,1)$, $H$ be a hypergraph, and $e\in E(H)$. If $\left\lceil p|e|\right\rceil=1$ then $e$ is \emph{highly flammable with respect to} $p$.
\end{definition}

As was the case with non-flammable edges, it is straightforward to determine whether or not an edge is highly flammable.

\begin{lemma}
An edge $e$ is highly flammable if and only if $|e|\leq\frac{1}{p}$.
\end{lemma}

\begin{proof}
Observe that $\left\lceil p|e|\right\rceil=1$ if and only if $p|e|\leq 1$, or equivalently, $|e|\leq \frac{1}{p}$.
\end{proof}

Singleton edges are the only edges that are both non-flammable and highly flammable. Although this is strange, it is not worthwhile to amend Definitions \ref{non_flam} and \ref{high_flam} since singleton edges have no effect on the game. In general, we disregard the existence of singleton edges for the sake of simplicity. For example, if we say ``all edges in $H$ are flammable,'' we really mean all non-singleton edges are flammable.

The remaining results in this section aim to compare the lazy and round-based versions of the game.

\begin{lemma}
Let $H$ be a connected hypergraph. Then there exists $p\in(0,1)$ such that $b_{L,p}(H)=1$. 
\end{lemma}

\begin{proof}
Choose $p$ small enough such that $p|e|\leq1$ for all $e\in E(H)$. Then every edge $e$ is highly flammable. Since $H$ is connected, a single vertex on fire in $H$ will cause the whole hypergraph to burn through subsequent propagation. So $b_{L,p}(H)=1$.
\end{proof}

\begin{lemma}
\label{one_implies_EH}
If $b_{L,p}(H)=1$ then $b_p(H)\leq |E(H)|+1$.
\end{lemma}

\begin{proof}
Let $\{v\}$ be a minimum lazy burning set for $H$ when using proportion $p$. If in any time step the fire does not propagate within an edge, then in fact $\{v\}$ is not a lazy burning set, so at least one edge catches fire every time step. Thus, it takes at most $|E(H)|$ time steps for $H$ to become fully burned. Now, the arsonist can burn $H$ in no more than $|E(H)|+1$ rounds by using the following strategy: burn $v$ as the first source (which takes one round), and then burn a redundant source every round for the rest of the game (which takes at most $|E(H)|$ rounds). This yields a burning sequence no longer than $|E(H)|+1$, so $b_p(H)\leq |E(H)|+1$.
\end{proof}

The bound in Lemma \ref{one_implies_EH} is tight -- consider the hypergraph consisting of a single edge containing all of its vertices. If the edge is highly flammable, then the round-based game takes two rounds. 

\begin{theorem}
\label{comparing_proportions}
Let $p,q\in(0,1)$ and $p\leq q$. Then for any hypergraph $H$, $b_{L,p}(H)\leq b_{L,q}(H)$ and $b_p(H)\leq b_q(H)$.
\end{theorem}

\begin{proof}
Observe that for any edge $e$, $\lceil p|e|\rceil\leq\lceil q|e|\rceil$. Hence, if at least $\lceil q|e|\rceil$ vertices are on fire in $e$, then it is also true that at least $\lceil p|e|\rceil$ vertices are on fire in $e$. Thus, if fire would propagate within $e$ when using proportion $q$, it would also propagate within $e$ when using proportion $p$. This means that any lazy burning set that is successful when using proportion $q$ is also successful when using proportion $p$ (if anything the propagation will be faster using proportion $p$). Therefore, $b_{L,p}(H)\leq b_{L,q}(H)$. 

Now, let $S$ be an optimal burning sequence for a hypergraph $H$ when using proportion $q$. Suppose the arsonist burns $H$ using proportion $p$ by following $S$. Recall that, when using proportion $p$, the propagation may happen faster than when using proportion $q$. Thus, when the arsonist tries to burn $H$ using proportion $p$ by following $S$, it is possible that they will be prompted to burn a vertex that was already on fire for at least one round, which is not allowed. When this occurs, the arsonist can simply deviate from $S$ for one round and choose any redundant source instead. Therefore, the arsonist can burn $H$ using proportion $p$ in no more than $|S|$ rounds, so $b_p(H)\leq |S|=b_q(H)$.
\end{proof}

\begin{theorem}
\label{lowerboundanyproportion}
Let $H$ be a hypergraph and $p\in(0,1)$. Then $\min\left\{\left\lceil p|e|\right\rceil \mid e\in E(H)\right\}\leq b_{L,p}(H)$. 
\end{theorem}

\begin{proof}
Clearly no fire can propagate if the lazy burning set contains $\min\left\{\left\lceil p|e|\right\rceil \mid e\in E(H)\right\}$ or fewer vertices.
\end{proof}

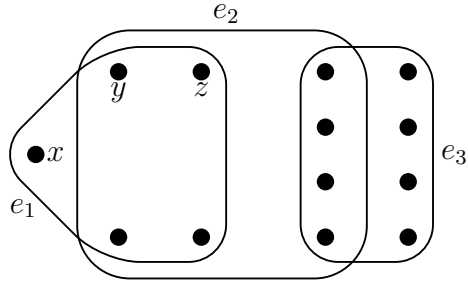
\begin{figure}[h]
\centering
\begin{tikzpicture}[scale=1.1]

\node (a) at (-1,1) {};
\fill [fill=black] (a) circle (0.105) node [below] {$z$};
\node (b) at (-2,1) {};
\fill [fill=black] (b) circle (0.105) node [below] {$y$};
\node (c) at (-3,0) {};
\fill [fill=black] (c) circle (0.105) node [right] {$x$};
\node (d) at (-1,-1) {};
\fill [fill=black] (d) circle (0.105) node [below] {};
\node (e) at (-2,-1) {};
\fill [fill=black] (e) circle (0.105) node [below] {};

\node (f) at (0.5,1) {};
\fill [fill=black] (f) circle (0.105) node [below] {};
\node (g) at (0.5,0.33) {};
\fill [fill=black] (g) circle (0.105) node [below] {};
\node (h) at (0.5,-0.33) {};
\fill [fill=black] (h) circle (0.105) node [below] {};
\node (i) at (0.5,-1) {};
\fill [fill=black] (i) circle (0.105) node [below] {};

\node (j) at (1.5,1) {};
\fill [fill=black] (j) circle (0.105) node [below] {};
\node (k) at (1.5,0.33) {};
\fill [fill=black] (k) circle (0.105) node [below] {};
\node (l) at (1.5,-0.33) {};
\fill [fill=black] (l) circle (0.105) node [below] {};
\node (m) at (1.5,-1) {};
\fill [fill=black] (m) circle (0.105) node [below] {};

\node () at (-3.15,-0.65) {$e_1$};
\node () at (-0.7,1.7) {$e_2$};
\node () at (2.06,0) {$e_3$};

\draw [line width=0.25mm] [black] [rounded corners=0.5cm] (-3.5,0)--(-2.2,1.3)--(-0.7,1.3)--(-0.7,-1.3)--(-2.2,-1.3)--cycle;
\draw [line width=0.25mm] [black] [rounded corners=0.5cm] (0.2,1.3)--(1.8,1.3)--(1.8,-1.3)--(0.2,-1.3)--cycle;
\draw [line width=0.25mm] [black] [rounded corners=0.7cm] (-2.5,1.5)--(1,1.5)--(1,-1.5)--(-2.5,-1.5)--cycle;

\end{tikzpicture}
\caption{An example that shows the bound in Theorem \ref{lowerboundanyproportion} is tight when $p\in (0,\frac{1}{2}]$.}
\label{halfburnlowerbound}
\end{figure}

The bound in Theorem \ref{lowerboundanyproportion} can be tight for any $p\in(0,1)$ -- simply consider a hypergraph consisting of a single edge that contains all of its vertices. A less trivial example that exhibits tightness for 
$p \in (0,\frac{1}{2}]$
can be seen in Figure \ref{halfburnlowerbound}. The hypergraph $H$ in Figure \ref{halfburnlowerbound} has 
$b_{L,p}(H) = \min\{\left\lceil p|e|\right\rceil\}$.
A minimum lazy burning set can be obtained by taking 
$\lceil 5p \rceil$
vertices of $e_1$.  



\begin{theorem}
\label{loose_ineq}
For any hypergraph $H$ and any proportion $p\in(0,1)$, $b_{L,p}(H)\leq b_p(H)$.
\end{theorem}

\begin{proof}
Fix any proportion $p$. If $(u_1,u_2,...,u_k)$ is a burning sequence for $H$, then $\{u_1,u_2,\ldots,u_k\}$ is a lazy burning set for $H$, so $b_{L,p}(H)\leq k$. In particular, we can choose an optimal burning sequence (i.e. one with $k=b_p(H)$), and then we have $b_{L,p}(H)\leq k=b_p(H)$.
\end{proof}

\begin{figure}[h]
\centering
\begin{tikzpicture}[scale=1.2]

\node (1) at (0,0) {};
\fill [fill=black] (1) circle (0.08) node [below] {};
\node (2) at (0,1) {};
\fill [fill=black] (2) circle (0.08) node [above] {};
\node (3) at (-1,-0.5) {};
\fill [fill=black] (3) circle (0.08) node [above right] {};
\node (4) at (-2,1) {};
\fill [fill=black] (4) circle (0.08) node [below right] {};
\node (5) at (-2,0) {};
\fill [fill=black] (5) circle (0.08) node [below] {};
\node (3) at (-1,1.5) {};
\fill [fill=black] (3) circle (0.08) node [above right] {};

\node (3) at (1.1,-0.5) {};
\fill [fill=black] (3) circle (0.08) node [above right] {};
\node (3) at (1.1,1.5) {};
\fill [fill=black] (3) circle (0.08) node [above right] {};
\node (3) at (2,0.5) {};
\fill [fill=black] (3) circle (0.08) node [above right] {};

\node (e1) at (-2.6,0.5) {$e_1$};
\node (e2) at (2.5,0.5) {$e_2$};

\draw [line width=0.25mm] [black] (-1,0.5) ellipse (1.4 cm and 1.25cm);
\draw [line width=0.25mm] [black] (0.95,0.5) ellipse (1.3 cm and 1.25cm);
\end{tikzpicture}

\caption{A hypergraph $H$ with one non-flammable edge and $b_{L,p}(H)=b_p(H)$ when $p=\frac{5}{6}$.}
\label{one_non_flam_edge}
\end{figure}
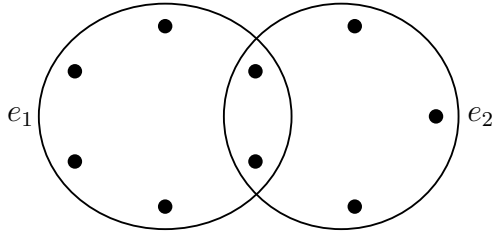

Note that if $H$ contains even one non-flammable edge, then we cannot say for certain that $b_{L,p}(H)$ is strictly less than $b_p(H)$. Observe the hypergraph $H$ in Figure \ref{one_non_flam_edge}. Suppose we burn $H$ using the proportion $p=\frac{5}{6}$. Then $e_1$ is flammable, and $e_2$ is non-flammable. Indeed, $b_{L,p}(H)=b_p(H)=8$.


\begin{theorem}
\label{tight_ineq}
Let $H$ be a hypergraph with no isolated vertices, and $p\in(0,1)$. If every edge in $H$ is flammable, then $b_{L,p}(H)< b_p(H)$.
\end{theorem}

\begin{proof}
Fix any proportion $p$, and let $(u_1,u_2,\ldots,u_{b_p(H)})$ be an optimal burning sequence for $H$ when burning using proportion $p$. It is clear that $\{u_1,u_2,\ldots, u_{b_p(H)}\}$ is a lazy burning set for $H$. We claim that $\{u_1,u_2,\ldots,u_{b_p(H)-1}\}$ is also a lazy burning set for $H$. Since $u_{b_p(H)}$ is the final source, burning $u_1,u_2,\ldots,u_{b_p(H)-1}$ one-by-one (as in the original game) or simultaneously (as a lazy burning set) will eventually result in all of $V(H)\setminus \{u_{b_p(H)}\}$ being burned. So, let us burn each vertex in $\{u_1,u_2,\ldots,u_{b_p(H)-1}\}$ simultaneously. All we need to show is that $u_{b_p(H)}$ will eventually burn through propagation. If $u_{b_p(H)}$ is a redundant source then clearly fire will propagate to $u_{b_p(H)}$. Otherwise, $u_{b_p(H)}$ is not a redundant source. But eventually all of $V(H)\setminus \{u_{b_p(H)}\}$ will burn through propagation, and $u_{b_p(H)}$ belongs to an edge since it is not isolated. All other vertices in the edge containing $u_{b_p(H)}$ are on fire, and thus $u_{b_p(H)}$ will catch on fire. Since $\{u_1,u_2,\ldots,u_{b_p(H)-1}\}$ is a lazy burning set for $H$ we have $b_{L,p}(H) \leq b_p(H)-1$, or equivalently, $b_{L,p}(H)<b_p(H)$.
\end{proof}

\begin{theorem}
\label{special_set_V}
Let $H$ be a hypergraph with no isolated vertices and $p\in(0,1)$. Let $V$ be a smallest set of vertices in $H$ such that, for every edge $e$ in $H$, at least $\lceil p|e|\rceil$ vertices in $e$ are also in $V$. Then $b_{L,p}(H)\leq |V|$ and $b_p(H)\leq |V|+1$.
\end{theorem}

\begin{proof}
First, clearly $V$ is a successful lazy burning set for $H$. After burning $V$ as a lazy burning set, every edge contains enough burned vertices to cause the fire to propagate. Since $H$ has no isolated vertices, every unburned vertex in $H$ catches fire in the first time step. Therefore, $b_{L,p}(H)\leq |V|$. 

In order to burn $H$ in no more than $ |V|+1$ rounds, the arsonist may use the following strategy. Each round, they burn an unburned vertex of $V$, so $V$ will be completely burned in at most $|V|$ rounds. In the following round, the rest of $H$ catches fire, and the arsonist burns a redundant source. This yields a burning sequence of length no greater than $|V|+1$, so $b_p(H)\leq |V|+1$.
\end{proof}


\begin{corollary}
\label{sum_round_based_upper_bound}
Let $H$ be a hypergraph with no isolated vertices, and let $p\in(0,1)$. Then $b_{L,p}(H)\leq \sum_{e\in E(H)}\lceil p|e|\rceil$ and $b_p(H)\leq 1+\sum_{e\in E(H)}\lceil p|e|\rceil$.
\end{corollary}

\begin{proof}
Let $V$ be a smallest set of vertices in $H$ such that, for every edge $e$ in $H$, at least $\lceil p|e|\rceil$ vertices in $e$ are also in $V$. Then clearly $|V|\leq \sum_{e\in E(H)}\lceil p|e|\rceil$, and the result follows by Theorem \ref{special_set_V}.
\end{proof}

\begin{corollary}
\label{at_most_k_round_based_upper_bound}
Let $H$ be a hypergraph with no isolated vertices whose largest edge contains $k$ vertices, and let $p\in(0,1)$. Then $b_{L,p}(H)\leq \lceil pk\rceil\cdot |E(H)|$ and $b_p(H)\leq \left(\lceil pk\rceil\cdot |E(H)|\right)+1$.
\end{corollary}

\begin{proof}
Keeping in mind that $|e|\leq k$ for each $e\in E(H)$, we calculate $\sum_{e\in E(H)}\lceil p|e|\rceil\leq \sum_{e\in E(H)}\lceil pk\rceil=\left(\lceil pk\rceil\cdot |E(H)|\right)$, and the result follows by Corollary \ref{sum_round_based_upper_bound}.
\end{proof}

\begin{corollary}
\label{k_unif_round_based_upper_bound}
Let $H$ be a $k$-uniform hypergraph and $p\in(0,1)$. Then $b_{L,p}(H)\leq \lceil pk\rceil\cdot |E(H)|$ and $b_p(H)\leq \left(\lceil pk\rceil\cdot |E(H)|\right)+1$.
\end{corollary}

The bounds in Theorem \ref{special_set_V}  and its corollaries (\ref{sum_round_based_upper_bound}, \ref{at_most_k_round_based_upper_bound}, and \ref{k_unif_round_based_upper_bound}) are all tight -- consider the $k$-uniform hypergraph in which no edges share a vertex. Now, we may combine the bounds in Theorem \ref{lowerboundanyproportion}, Theorem \ref{tight_ineq}, and Corollary \ref{sum_round_based_upper_bound} to get the following result.

\begin{theorem}
\label{combined_bounds}
Let $H$ be a hypergraph with no isolated vertices, and let $p\in(0,1)$ such that every edge in $H$ is flammable. Then $$\min\left\{\left\lceil p|e|\right\rceil \mid e\in E(H)\right\}\leq b_{L,p}(H)< b_p(H)\leq 1+\sum_{e\in E(H)}\lceil p|e|\rceil.$$ 
\end{theorem}

The bounds in Theorem \ref{combined_bounds} can all be tight simultaneously -- consider the hypergraph consisting of a single edge containing all of its vertices.


Our next major result is a tight upper bound on the lazy burning number of a hypergraph when using a proportion of the form $p=\frac{1}{n}$; see Theorem \ref{pproptheorem}. The following elementary result is used in its proof.

\begin{lemma}
\label{countinglemma}
If $a$ and $b$ are real numbers, then $\lceil a\rceil + \lfloor b\rfloor \leq \lceil a+b\rceil$.
\end{lemma}

\begin{theorem}
\label{pproptheorem}
Let $H$ be a connected hypergraph and $n\in\mathbb{N}$. Then $b_{L,\frac{1}{n}}(H)\leq \left\lceil\frac{|V(H)|}{n}\right\rceil$. 
\end{theorem}

\begin{proof}
We use induction on $|E(H)|$, assuming singleton edges are not allowed. The bound in this theorem will still be valid for hypergraphs containing singleton edges, as they have no effect on the burning game.

\noindent \emph{Base Case.} If $|E(H)|=1$ then $H$ consists of a single edge $e$ containing all of $V(H)$. Thus, $b_{L,\frac{1}{n}}(H)=\left\lceil\frac{|e|}{n}\right\rceil=\left\lceil\frac{|V(H)|}{n}\right\rceil$.

\noindent \emph{Inductive Hypothesis.} Suppose that for some $k\in\mathbb{N}$, every connected hypergraph $H$ with $|E(H)|=k$ satisfies the bound $b_{L,\frac{1}{n}}(H)\leq \left\lceil\frac{|V(H)|}{n}\right\rceil$.

\noindent \emph{Inductive Step.} Let $G$ be a connected hypergraph with $k+1$ edges. Observe that, since $G$ is connected, there must exist some edge $e\in E(G)$ such that $(V(G),E(G)\setminus\{e\})$ is connected, except for possibly some isolated vertices (degree $1$ vertices in $e$ would become isolated vertices in $(V(G),E(G)\setminus\{e\})$). Choose any such $e\in E(G)$, and write $e=\{u_1,\ldots,u_x\}\cup\{v_1,\ldots,v_y\}$ where each $u_i$ is a vertex in $e$ of degree greater than $1$, and each $v_j$ is a degree $1$ vertex in $e$. Since $G$ is connected, $e$ must contain at least one vertex of degree greater than $1$, so $x\geq 1$. Also, clearly $y\geq 0$. Now, define $H$ as the hypergraph with $V(H)=V(G)\setminus\{v_1,\ldots,v_y\}$ and $E(H)=E(G)\setminus\{e\}$. Then $H$ is a connected hypergraph with $k$ edges, so $b_{L,\frac{1}{n}}(H)\leq \left\lceil\frac{|V(H)|}{n}\right\rceil = \left\lceil\frac{|V(G)|-y}{n}\right\rceil$. Let $S$ be a minimum lazy burning set in $H$, so $|S|\leq\left\lceil\frac{|V(G)|-y}{n}\right\rceil$.

Now, we claim $S^\prime=S\cup\{v_1,\ldots,v_{\left\lfloor \frac{y}{n}\right\rfloor}\}$ is a lazy burning set for $G$. All of $V(G)\setminus\{v_1,\ldots,v_y\}$ will catch fire because $S\subseteq S^\prime$. In particular, $u_1,\ldots,u_x,v_1,\ldots,$ and $v_{\left\lfloor \frac{y}{n}\right\rfloor}$ will all eventually be on fire either due to propagation or because they were included in $S^\prime$. We now show that the rest of $e$ will catch fire. In particular, we must show that $x+\left\lfloor \frac{y}{n}\right\rfloor \geq \left\lceil\frac{|e|}{n}\right\rceil$, since $x+\left\lfloor \frac{y}{n}\right\rfloor$ vertices are sure to be on fire in $e$. In our calculations, we will use the well-known identities $\left\lfloor\frac{a}{b}\right\rfloor=\left\lceil\frac{a+1}{b}\right\rceil-1$ and $\lceil c\rceil +a=\lceil c+a\rceil$ for any $a,b\in\mathbb{N}$ and $c\in\mathbb{R}$.

First, observe $1\leq x\implies (n-1)\leq x(n-1)\implies 0\leq x(n-1)+1-n\implies 0\leq\frac{(n-1)x+1-n}{n}$, which will be used in the final line of the calculations below. Now, \begin{align*}
x+\left\lfloor \frac{y}{n}\right\rfloor &= x+\left\lceil \frac{y+1}{n}\right\rceil-1 \\
&=\left\lceil x+\frac{y+1}{n}-1\right\rceil \\
&=\left\lceil \frac{nx+y+1-n}{n}\right\rceil \\
&=\left\lceil \frac{x+y}{n}+\frac{(n-1)x+1-n}{n}\right\rceil \\
&=\left\lceil \frac{|e|}{n}+\frac{(n-1)x+1-n}{n}\right\rceil \\
&\geq \left\lceil \frac{|e|}{n}\right\rceil.
\end{align*}

Therefore $e$ will cause fire to spread to each of the vertices $v_{\left\lfloor \frac{y}{n} \right\rfloor+1},\ldots,v_y$. This will leave $G$ completely burned, so $S^\prime$ is a lazy burning set for $G$. Now, using Lemma \ref{countinglemma}, we calculate $$b_{L,\frac{1}{n}}(G) \leq |S^\prime|=|S|+\left\lfloor\frac{y}{n}\right\rfloor \leq \left\lceil\frac{|V(G)|-y}{n}\right\rceil+\left\lfloor\frac{y}{n}\right\rfloor \leq \left\lceil\frac{|V(G)|-y}{n}+\frac{y}{n}\right\rceil=\left\lceil\frac{|V(G)|}{n}\right\rceil.$$
\end{proof}

\begin{figure}[h]
\centering
\begin{tikzpicture}[scale=1.1]

\node (a) at (-1,0) {};
\fill [fill=black] (a) circle (0.105) node [below] {};
\node (b) at (-2,1) {};
\fill [fill=black] (b) circle (0.105) node [below] {$q$};
\node (c) at (-3,1) {};
\fill [fill=black] (c) circle (0.105) node [below] {$p$};
\node (d) at (-3,-1) {};
\fill [fill=black] (d) circle (0.105) node [above] {$r$};
\node (e) at (-2,-1) {};
\fill [fill=black] (e) circle (0.105) node [below] {};

\node (f) at (1,0) {};
\fill [fill=black] (f) circle (0.105) node [below] {};
\node (g) at (2,1) {};
\fill [fill=black] (g) circle (0.105) node [below] {$s$};
\node (h) at (3,1) {};
\fill [fill=black] (h) circle (0.105) node [below] {$t$};
\node (i) at (3,-1) {};
\fill [fill=black] (i) circle (0.105) node [below] {};

\node (j) at (2,-1) {};
\fill [fill=black] (j) circle (0.105) node [below] {};

\node () at (-3.53,0) {$e_1$};
\node () at (0,0.72) {$e_2$};
\node () at (3.58,0) {$e_3$};

\draw [line width=0.25mm] [black] [rounded corners=0.5cm] (0.6,0)--(1.8,1.3)--(3.3,1.3)--(3.3,-1.3)--(1.8,-1.3)--cycle;
\draw [line width=0.25mm] [black] [rounded corners=0.5cm] (-0.6,0)--(-1.8,1.3)--(-3.3,1.3)--(-3.3,-1.3)--(-1.8,-1.3)--cycle;
\draw [line width=0.25mm, black ] (0,0) ellipse (1.5cm and 0.5cm);

\end{tikzpicture}
\caption{An example that shows the bound in Theorem \ref{pproptheorem} can be tight.}
\label{halfpropupperbound}
\end{figure}
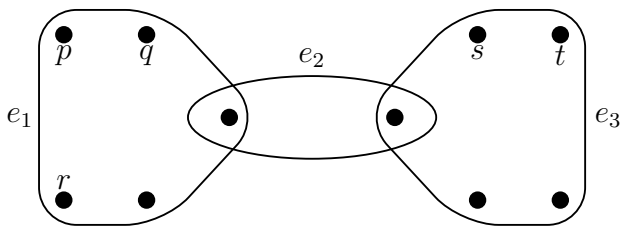

The bound in Theorem \ref{pproptheorem} can be tight for any $n\in\mathbb{N}\setminus\{1\}$ -- consider a connected hypergraph consisting of one edge that contains all of its vertices. A less trivial example that exhibits tightness can be seen in Figure \ref{halfpropupperbound}. The hypergraph $H$ in Figure \ref{halfpropupperbound} has $b_{L,\frac{1}{2}}(H)=5=\left\lceil\frac{|V(H)|}{2}\right\rceil$, as $\{p,q,r,s,t\}$ is a minimum lazy burning set when $p=\frac{1}{2}$. This hypergraph also has $b_{L,\frac{1}{4}}(H)=3=\left\lceil\frac{|V(H)|}{4}\right\rceil$, as $\{p,q,s\}$ is a minimum lazy burning set when $p=\frac{1}{4}$. This yields an example where the bound in Theorem \ref{pproptheorem} is tight for $p=\frac{1}{n}$ with $n>2$.


\begin{conjecture}
\label{ppropconj1}
Let $H$ be a connected hypergraph and $p\in(0,1)\cap\mathbb{Q}$ such that $H$ contains no non-flammable edges. Then $b_{L,p}(H)\leq \left\lceil p|V(H)|\right\rceil$. 
\end{conjecture}

Note that Conjecture \ref{ppropconj1} would imply the following Conjecture \ref{ppropconj2}. If we are given $p\in (0,1)\setminus\mathbb{Q}$, then we can simply choose $p^\prime\in (0,1)\cap\mathbb{Q}$ close enough to $p$ such that $b_{L,p}(H)=b_{L,p^\prime}(H)$ and $\left\lceil p|V(H)|\right\rceil=\left\lceil p^\prime |V(H)|\right\rceil$, and then apply the bound from Conjecture \ref{ppropconj1}.

\begin{conjecture}
\label{ppropconj2}
Let $H$ be a connected hypergraph and $p\in(0,1)$ such that $H$ contains no non-flammable edges. Then $b_{L,p}(H)\leq \left\lceil p|V(H)|\right\rceil$. 
\end{conjecture}

In order to use the same proof technique as in Theorem \ref{pproptheorem} to prove Conjecture \ref{ppropconj1}, we would need to prove the following for at least one edge $e\in E(G)$ (where $G$ is from the induction step): If $p\in(0,1)\cap\mathbb{Q}$, $x\in\mathbb{N}$, $y\in\mathbb{N}\cup\{0\}$, and $|e|=x+y\geq\frac{1}{1-p}$ ($e$ has $x$ vertices of degree greater than $1$, $y$ degree $1$ vertices, and is flammable), then $x+\lfloor py\rfloor \geq \lceil p|e|\rceil$. This is required for the step where we show that $S^\prime=S\cup\{ v_1,\ldots,v_{\lfloor py\rfloor}\}$ eventually causes the vertices $v_{\lfloor py\rfloor+1},\ldots,v_y$ to catch fire, and hence it is a lazy burning set for $G$. However, counterexamples exist when $p$ is not of the form $\frac{1}{n}$. For instance, let $p=\frac{13}{20}$ (notably above $\frac{1}{2}$), and suppose $e$ has $x=2$ vertices of degree greater than one and $y=3$ degree-one vertices. Then $\frac{1}{1-p}=\frac{1}{1-\frac{13}{20}}=\frac{20}{7}<5=|e|$, so $e$ is flammable. Also, $\lceil p|e|\rceil= \left\lceil \frac{13}{20}\cdot 5\right\rceil=4$, so four vertices are required to be on fire for fire to propagate within $e$. But $x+\lfloor py\rfloor=2+\left\lfloor\frac{13}{20}\cdot3\right\rfloor=2+\left\lfloor\frac{39}{20}\right\rfloor=2+1=3$, so $S^\prime$ will not cause fire to propagate within $e$. Another counterexample arises when $p=\frac{2}{5}$ (notably below $\frac{1}{2}$), $x=1$, and $y=2$. 


By combining Theorems \ref{comparing_proportions} and \ref{pproptheorem}, we can get the following upper bound on $b_{L,p}(H)$ for $p\in\left(0,\frac{1}{2}\right]$ not necessarily of the form $\frac{1}{n}$. 

\begin{corollary}
Let $p\in\left(0,\frac{1}{2}\right]$, and let $n$ be the largest number in $\mathbb{N}\setminus\{1\}$ such that $p\leq \frac{1}{n}$. Then $b_{L,p}(H)\leq \left\lceil\frac{|V(H)|}{n}\right\rceil$.
\end{corollary}

\begin{theorem}
\label{comparing_to_original}
If every edge in $H$ is flammable when burning with proportion $p$ then $b_{L,p}(H)\leq b_L(H)$ and $b_p(H)\leq b(H)$.
\end{theorem}

\begin{proof}
One may use the same argument as in the proof of Theorem \ref{comparing_proportions}, since the fire propagates more efficiently within every edge in the proportion-based variant of the game when all edges are flammable.
\end{proof}

The differences $b(H)-b_p(H)$ and $b_L(H)-b_{L,p}(H)$ can be arbitrarily large --  just consider a hypergraph $H$ with one edge that contains all of its vertices.


By combining Lemma \ref{non_flam_iff} and Theorem \ref{comparing_to_original} we get the following result.

\begin{corollary}
Let $\ell$ be the size of a smallest non-singleton edge in $H$. If $\ell\geq\frac{1}{1-p}$ then $b_{L,p}(H)\leq b_L(H)$ and $b_p(H)\leq b(H)$.
\end{corollary}

In the original model hypergraph burning, the values of $b(H)$ and $b_L(H)$ can differ by an arbitrarily large amount. More specifically, for any $\ell\in\mathbb{N}$, there exists a hypergraph $H$ such that $b(H)-b_L(H)>\ell$ and $\frac{b(H)}{b_L(H)}>\ell$ \cite{mythesis}. Although it seems $b_{L,p}(H)$ and $b_p(H)$ get closer together as $p$ approaches one (see Remark \ref{all_non_flam_cor}), an analogous result does indeed exist for proportion-based hypergraph burning. We first need the following definition.

\begin{definition}
If $V(H)=\{v_1,\ldots,v_n\}$ and $E(H)=\big{\{} \{v_1,\dots,v_k\},\{v_2,\ldots,v_{k+1}\},$\\ $\ldots,\{v_{n-k+1},\ldots,v_n\} \big{\}}$ then $H$ is a $k$-uniform \emph{tight path}. 
\end{definition}

\begin{theorem}
Let $p\in(0,1)$ and $\ell\in\mathbb{N}$. There exists a hypergraph $H$ such that $b_p(H)-b_{L,p}(H)>\ell$ and $\frac{b_p(H)}{b_{L,p}(H)}>\ell$.
\end{theorem}

\begin{proof}
Choose $k\in\mathbb{N}$ large enough so that $\lceil pk\rceil<k$, and consider the family of $k$-uniform tight paths. All such hypergraphs have lazy burning number $\lceil pk\rceil$ when burning with proportion $p$. Indeed, one may choose $\lceil pk\rceil$ vertices in an edge $e$ as the lazy burning set, so that in the first time step, $e$ catches completely on fire. Now, consider the two edges ``adjacent'' to $e$ (i.e. the two unique edges that share exactly $k-1$ vertices with $e$). Since $\lceil pk\rceil\leq k-1$, both of these edges will catch completely on fire in the second time step. Clearly, the propagation will continue and eventually leave the hypergraph completely burned.

However, the burning number (using proportion $p$) monotonically increases with the number of edges in this family of hypergraphs. The burning number for this family of hypergraphs is unbounded, and the lazy burning number is a fixed, finite number, so the result follows.
\end{proof}

Observe that even if $p$ is low, it is not necessarily the case
that $b_{L,p}(H)$ and $b_p(H)$ differ by a large amount. For example, one may consider a hypergraph consisting of a single edge containing all of its vertices. For very low $p$, one may take any such hypergraph with two or more vertices in which the edge is highly flammable. Then $b_{L,p}(H)=1$ and $b_p(H)=2$.

We close this section with a brief discussion of complexity. It is easy to come up with a polynomial-time algorithm which takes a proportion $p\in(0,1)$ and a sequence of vertices $S$ in a hypergraph $H$ as input, and determines if $S$ is a valid burning sequence for $H$ when using proportion $p$. There is a similar algorithm for the lazy burning game which determines if a given subset of $V(H)$ is a lazy burning set (when burning using a given proportion $p$) in polynomial time. Hence, both games are in {NP}. Graph burning is {NP}-complete \cite{burning_is_hard}, and it is a special case of proportion-based hypergraph burning. This is because burning a $2$-uniform hypergraph with any proportion $p\in(0,\frac{1}{2}]$ is equivalent to graph burning. Thus, proportion-based hypergraph burning is also {NP}-complete. The question of whether proportion-based lazy hypergraph burning is {NP}-complete is still open. 

\section{Results on the Burning Distribution}

Since hypergraphs are discrete structures, being able to choose any real proportion $p\in(0,1)$ seems like ``overkill'' in the sense that it is a redundant amount of accuracy. Given a hypergraph $H$, there are uncountably many $p$ that all yield the same value for $b_p(H)$, since for any $p\in(0,1)$, $b_p(H)\in\{1,2,\ldots,|V(H)|\}$ (and similarly for $b_{L,p}(H)$). Moreover, in light of Theorem \ref{comparing_proportions}, there must surely be ``cut-off points'' for $p$. That is, there must be proportions $q_k\in(0,1)$ such that if $p<q_k$ then $b_p(H)<k$, and if $p>q_k$ then $b_p(H)\geq k$ (and similarly for the lazy version). We therefore introduce the following definition, which gives insight into the range of possible values for $b_p(H)$ and $b_{L,p}(H)$ given a hypergraph $H$ (where $p$ ranges over $(0,1)$), as well as what the aforementioned ``cut-off points'' are for $p$ when burning $H$.

\begin{definition}
Let $H$ be a hypergraph with $|V(H)|=n$. The \emph{burning distribution} of $H$ is a partition $\mathcal{P}$ of $(0,1)$ into $n$ (possibly empty) intervals $P_1,P_2,\ldots,P_n$ such that $b_p(H)=k$ if and only if $p\in P_k$. The \emph{lazy burning distribution} is defined analagously, and denoted $\mathcal{P}_L$. As a convention, we will denote the intervals in $\mathcal{P}_L$ by $Q_1,Q_2,\ldots,Q_n$.
\end{definition}


\begin{lemma}
For any $n\in\mathbb{N}$, there is a hypergraph $H$ with $|V(H)|=n$ such that every interval in the lazy burning distribution is nonempty.
\end{lemma}

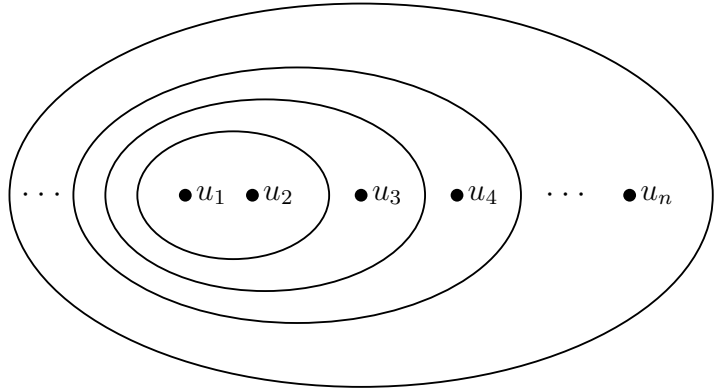
\begin{figure}
\centering
\begin{tikzpicture}[scale=0.85]
\draw [line width=0.25mm, black ] (-1.5,0) ellipse (1.5cm and 1cm);
\draw [line width=0.25mm, black ] (-1,0) ellipse (2.5cm and 1.5cm);
\draw [line width=0.25mm, black ] (-0.5,0) ellipse (3.5cm and 2cm);
\draw [line width=0.25mm, black ] (0.5,0) ellipse (5.5cm and 3cm);
\node at (-4.5,0) {$\cdot\cdot\cdot$};   
\node at (3.7,0) {$\cdot\cdot\cdot$};   
\node (u1) at (-2.25,0) {};
\fill [fill=black] (u1) circle (0.1) node [right] {$u_1$};
\node (u2) at (-1.2,0) {};
\fill [fill=black] (u2) circle (0.1) node [right] {$u_2$};   
\node (u3) at (0.5,0) {};
\fill [fill=black] (u3) circle (0.1) node [right] {$u_3$};    
\node (u4) at (2,0) {};
\fill [fill=black] (u4) circle (0.1) node [right] {$u_4$};    
\node (ul) at (4.7,0) {};
\fill [fill=black] (ul) circle (0.1) node [right] {$u_n$}; 
\end{tikzpicture}
\caption{A hypergraph in which every interval in the lazy burning distribution is nonempty.}
\label{lazy_dist_all}
\end{figure}

\begin{proof}
Consider the hypergraph $H$ seen in Figure \ref{lazy_dist_all}, which has $V(H)=\{u_1,\ldots,u_n\}$ and $E(H)=\big{\{}\{u_1,u_2\},\{u_1,u_2,u_3\},\{u_1,u_2,u_3,u_4\},\ldots,\{u_1,u_2,\ldots,u_n\}\big{\}}$. There exists a lazy burning set of size $1$ if and only if $p\in(0,\frac{1}{2}]$, so $Q_1=(0,\frac{1}{2}]$ (one example of such a lazy burning set is $\{u_1\}$). There exists a minimum lazy burning set of size $2$ if and only if $p\in(\frac{1}{2},\frac{2}{3}]$, so $Q_2=(\frac{1}{2},\frac{2}{3}]$ (one example of such a lazy burning set is $\{u_1,u_2\}$). In general, for any $k<n$, $\{u_1,u_2,\ldots,u_k\}$ is a minimum lazy burning set if and only if $p\in(\frac{k-1}{k},\frac{k}{k+1}]$, so $Q_k=(\frac{k-1}{k},\frac{k}{k+1}]$. Finally, a minimum lazy burning set of size $n$ exists if and only if $p>\frac{n-1}{n}$, since then every edge is non-flammable. Thus, $Q_n=(\frac{n-1}{n},1)$.
\end{proof}

The following remark essentially says that, if we increase the proportion from $p_1$ to $p_2$ such that some edge $e\in E(H)$ burns ``less easily'' using proportion $p_2$, it may still be true that $b_{L,p_1}(H)= b_{L,p_2}(H)$ and $b_{p_1}(H)= b_{p_2}(H)$. That is, the (lazy) burning number does not change even though the fire propagates differently within some edge.

\begin{remark}
There exists a hypergraph $H$ and proportions $p_1,p_2\in(0,1)$ such that $\lceil p_1|e|\rceil<\lceil p_2|e|\rceil$ for some edge $e\in E(H)$, but $b_{L,p_1}(H)= b_{L,p_2}(H)$ and $b_{p_1}(H)= b_{p_2}(H)$.
\end{remark}

\begin{figure}[h]
\centering
\begin{tikzpicture}[scale=1.2]

\node (1) at (0,0) {};
\fill [fill=black] (1) circle (0.08) node [below] {};
\node (2) at (-1,1) {};
\fill [fill=black] (2) circle (0.08) node [above] {};
\node (3) at (-1,-1) {};
\fill [fill=black] (3) circle (0.08) node [above right] {};
\node (4) at (-2.2,0.6) {};
\fill [fill=black] (4) circle (0.08) node [below right] {};
\node (5) at (-2.2,-0.6) {};
\fill [fill=black] (5) circle (0.08) node [below] {};
\node (3) at (2,0) {};
\fill [fill=black] (3) circle (0.08) node [above right] {};
\node (3) at (1.3,1) {};
\fill [fill=black] (3) circle (0.08) node [above right] {};
\node (3) at (1.3,-1) {};
\fill [fill=black] (3) circle (0.08) node [above right] {};

\node (e1) at (0.5,0) {$e_1$};
\node (e2) at (2.7,0) {$e_2$};

\draw [line width=0.25mm] [black] (-1.2,0) ellipse (1.5 cm and 1.25cm);
\draw [line width=0.25mm] [black] (-0.3,0) ellipse (2.8 cm and 1.7cm);
\end{tikzpicture}
\caption{An example where increasing the proportion such that an edge behaves differently does not necessarily increase the (lazy) burning number.}
\label{one_edge_changes}
\end{figure}
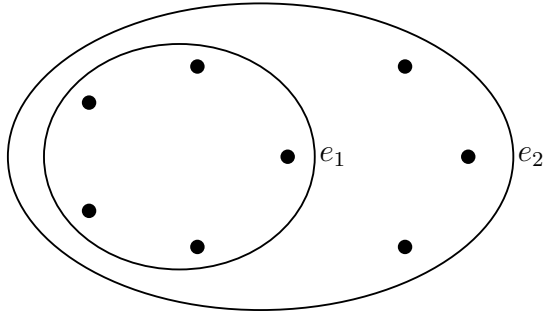

\begin{proof}
Observe the hypergraph $H$ in Figure \ref{one_edge_changes}, and set $p_1=\frac{3}{10}$ and $p_2=\frac{2}{5}$. We calculate $\lceil p_1 |e_1|\rceil=\left\lceil\frac{3}{10}\cdot 5\right\rceil=2$, $\lceil p_2 |e_1|\rceil=\left\lceil\frac{2}{5}\cdot 5\right\rceil=2$, $\lceil p_1 |e_2|\rceil=\left\lceil\frac{3}{10}\cdot 8\right\rceil=3$, and $\lceil p_2 |e_2|\rceil=\left\lceil\frac{2}{5}\cdot 8\right\rceil=4$. Whether we burn with proportion $p_1$ or $p_2$, edge $e_1$ behaves the same. However, when we burn using proportion $p_1$, three vertices on fire in $e_2$ cause fire to propagate, and when we burn using proportion $p_2$, four vertices on fire in $e_2$ cause fire to propagate. So edge $e_2$ behaves differently depending on which proportion we use. But $b_{L,p_1}(H)=b_{L,p_2}(H)=2$ and $b_{p_1}(H)=b_{p_2}(H)=4$.
\end{proof}


Next, we investigate the nature of the (lazy) burning distribution on $k$-uniform hypergraphs.

\begin{lemma}
\label{k_unif_Q1_Qn}
In a connected $k$-uniform hypergraph $H$ on $n$ vertices, $Q_1=\left(0,\frac{1}{k}\right]$ and $Q_n=\left(1-\frac{1}{k},1\right)$.
\end{lemma}

\begin{proof}
First, to find $Q_1$ observe that $$p\in \left(0,\frac{1}{k}\right] \iff p\leq\frac{1}{k}\iff k\leq\frac{1}{p}\iff \textnormal{every edge in}\ H\ \textnormal{is highly flammable}.$$ Indeed, since $H$ is $k$-uniform and connected, $b_{L,p}(H)=1$ exactly when every edge in $H$ is highly flammable, so $Q_1=\left(0,\frac{1}{k}\right]$. Now, to find $Q_n$ observe that $$p\in \left(1-\frac{1}{k},1\right)\iff p>1-\frac{1}{k}\iff k<\frac{1}{1-p}\iff \textnormal{every edge in}\ H\ \textnormal{is non-flammable}.$$ Indeed, $b_{L,p}(H)=n$ exactly when every edge in $H$ is non-flammable, so $Q_n=\left(1-\frac{1}{k},1\right)$.
\end{proof}

\begin{lemma}
\label{k_unif_Pn}
If $H$ is a simple $k$-uniform hypergraph with $n$ vertices and more than one edge, then $P_1=\emptyset$ and $P_n=\left(1-\frac{1}{k},1\right)$.
\end{lemma}

\begin{proof}
First, any hypergraph with two or more vertices will take at least two rounds to burn, regardless of the proportion $p$. Hence, $P_1=\emptyset$.

Clearly $\left(1-\frac{1}{k},1\right)\subseteq P_n$, since if $p\in \left(1-\frac{1}{k},1\right)$ then every edge in $H$ is non-flammable and hence $b_p(H)=n$. 

Now, suppose $p\leq1-\frac{1}{k}$. Then $k\geq\frac{1}{1-p}$, so every edge in $H$ is flammable. The arsonist can begin the game by burning $\lceil pk\rceil$ vertices in an edge $e$, and then burning vertices in $V(H)\setminus e$ for the rest of the game. Thus, in round $\lceil pk\rceil+1$, at least one vertex catches fire through propagation in $e$, and the arsonist chooses a non-redundant source. There is a burning sequence which avoids burning every vertex in $H$ as a source, so $b_p(H)<n$ and thus $p\notin P_n$. We have proven that if $p\notin\left(1-\frac{1}{k},1\right)$ then $p\notin P_n$, or equivalently, $P_n\subseteq\left(1-\frac{1}{k},1\right)$.
\end{proof}

\begin{remark}
Let $H$ be a single edge containing $k>1$ vertices. Then $Q_1=\left(0,\frac{1}{k}\right]$, $Q_2=\left(\frac{1}{k},\frac{2}{k}\right]$, $\ldots$ ,$Q_{k-1}=\left(\frac{k-2}{k},\frac{k-1}{k}\right]$, $Q_{k}=\left(\frac{k-1}{k},1\right)$, and $P_1=\emptyset$, $P_2=\left(0,\frac{1}{k}\right]$, $P_3=\left(\frac{1}{k} ,\frac{2}{k} \right]$, $\ldots$ ,$P_{k-1}=\left(\frac{k-3}{k},\frac{k-2}{k}\right]$, $P_k=\left(\frac{k-2}{k},1\right)$.
\end{remark}

\begin{proof}
Clearly for any $1\leq i<k$, $b_{L,p}(H)=i$ if and only if $p\in\left(\frac{i-1}{k},\frac{i}{k}\right]$, since then $i$ vertices on fire in the edge cause fire to propagate. The only case when $b_{L,p}(H)=k$ is when the edge is non-flammable, which occurs exactly when $p>\frac{k-1}{k}$.

Now, the optimal strategy in the round-based game is to burn enough vertices to cause fire to propagate, and then spend one more round burning a redundant source. Hence, for any $1\leq i<k-1$, the game takes $i+1$ rounds if and only if $i$ vertices on fire cause the fire to propagate. Thus, for  $1\leq i<k-1$, $P_{i+1}=\left(\frac{i-1}{k},\frac{i}{k}\right]$. Observe that $P_1$ is empty, since even if the edge is highly flammable, the game still takes two rounds. Finally, observe that the game takes $k$ rounds if and only if the edge is non-flammable or $k-1$ vertices cause the fire to propagate. Therefore, $P_k=\left(\frac{k-2}{k},\frac{k-1}{k}\right]\cup \left(\frac{k-1}{k},1\right)$.
\end{proof}

The following theorem uses a common ``shadow strategy'' technique in its proof. To avoid confusion, we denote the burning game on $H$ using the burning sequence $S$ and proportion $p$ by $BG_p(H,S)$. Similarly, we denote the lazy burning game on $H$ with lazy burning set $T\subseteq V(H)$ and proportion $p$ by $LBG_p(H,T)$. 

\begin{theorem}
\label{k_unif_cutoff}
Let $H$ be a $k$-uniform hypergraph with $n$ vertices and multiple distinct edges. Also let $m\in\{1,2,\ldots,k-1\}$, $p_1=\frac{m}{k}$, and $p_2\in\left(\frac{m}{k},\frac{m+1}{k}\right)$. Then $b_{L,p_1}(H)<b_{L,p_2}(H)$ and $b_{p_1}(H)<b_{p_2}(H)$. 
\end{theorem}

\begin{proof}
First, observe that exactly $m$ vertices on fire in any edge will cause fire to propagate when burning $H$ with proportion $p_1$, and exactly $m+1$ vertices on fire in any edge will cause fire to propagate when burning $H$ with proportion $p_2$.

Let $L$ be a minimum lazy burning set for $H$ when burning with proportion $p_2$. Consider any edge $e$ that fully catches fire in the first time step. Since $L$ is of minimum size, exactly $m+1$ vertices of $e$ are in the lazy burning set $L$. Let $v$ be any vertex in $e\cap L$. We claim $L^\prime=L\setminus\{v\}$ is a lazy burning set for $H$ when burning with proportion $p_1$. Indeed, $e$ will fully catch fire in the first propagation step in $LBG_{p_1}(H,L^\prime)$ since $m$ vertices are set on fire in $e$ as part of the lazy burning set $L^\prime$. Now, since $p_1<p_2$, if fire would propagate in an edge of $H$ using proportion $p_2$, then it would also propagate using proportion $p_1$. Hence, the propagation that occurs in the second time step or later in $LBG_{p_2}(H,L)$ will also occur in $LBG_{p_1}(H,L^\prime)$. So $L^\prime$ is a lazy burning set for $H$ when using proportion $p_1$. Thus, $b_{L,p_1}(H)\leq |L^\prime|<|L|=b_{L,p_2}(H)$.

Now, let $S=(u_1,u_2,\ldots,u_\ell)$ be an optimal burning sequence for $H$ when using proportion $p_2$. Let $u_q$ be the lowest-indexed source that is the $(m+1)^{th}$ vertex to be burned in some edge $e$. Thus, in $BG_{p_2}(H,S)$, $e$ is the first edge to catch completely on fire, and this occurs in round $q+1$. Let $S^\prime$ be the sequence of sources in $H$ (using proportion $p_1$) that is constructed in the following way: for each round $i<q$, burn $u_i$ as a source in $S^\prime$. Then, for each round $i\geq q$, burn $u_{i+1}$. If in any round the desired source is already on fire, simply burn a redundant source that round instead. Informally, we are following $S$ as closely as possible, but ``skipping'' $u_q$, and stopping the process if $H$ becomes fully burned at any point. We will show that $H$ is fully burned in $BG_{p_1}(H,S^\prime)$ no later than round $\ell-1$.

At the end of rounds $1,2,\ldots,q-1$, the vertices on fire in $BG_{p_2}(H,S)$ are a subset of those on fire in $BG_{p_1}(H,S^\prime)$. At the end of round $q-1$ in $BG_{p_2}(H,S)$, there are $m$ vertices on fire in the edge $e$. Thus, at the end of round $q-1$ in $BG_{p_1}(H,S^\prime)$, there are at least $m$ vertices on fire in the edge $e$. So $e$ catches fully on fire in $BG_{p_1}(H,S^\prime)$ in round $q$ (or possibly earlier). Therefore, the vertices on fire in $BG_{p_2}(H,S)$ at the end of round $q+1$ are a subset of those on fire in $BG_{p_1}(H,S^\prime)$ at the end of round $q$. Indeed, since the fire propagates more efficiently in $BG_{p_1}(H,S^\prime)$ than in $BG_{p_2}(H,S)$, the set of vertices on fire at the end of round $i$ in $BG_{p_2}(H,S)$ will be a subset of those on fire at the end of round $i-1$ in $BG_{p_1}(H,S^\prime)$ for the remainder of the game. Since $H$ is fully burned at the end of round $\ell$ in $BG_{p_2}(H,S)$, $H$ must therefore be fully burned in $BG_{p_1}(H,S^\prime)$ at the end of round $\ell-1$ (or possibly earlier). Hence, $b_{p_1}(H)\leq \ell-1<\ell=b_{p_2}(H)$.
\end{proof}

\begin{corollary}
\label{k_unif_cutoff_cor}
Let $H$ be a $k$-uniform hypergraph with $n$ vertices and multiple distinct edges. Also let $m\in\{1,2,\ldots,k-1\}$ and $p_1\leq\frac{m}{k}<p_2$. Then $b_{L,p_1}(H)<b_{L,p_2}(H)$ and $b_{p_1}(H)<b_{p_2}(H)$. 
\end{corollary}

\begin{proof}
Let $p_3\in \left(\frac{m}{k},\frac{m+1}{k}\right)$ such that $p_3\leq p_2$. By Theorem \ref{comparing_proportions} and Theorem \ref{k_unif_cutoff}, we have $b_{L,p_1}(H)\leq b_{L,\frac{m}{k}}< b_{L,p_3}\leq b_{L,p_2}$, and $b_{p_1}(H)\leq b_{\frac{m}{k}}< b_{p_3}\leq b_{p_2}$.
\end{proof}

We are finally ready to determine the nature of the (lazy) burning distribution in a $k$-uniform hypergraph, which can be seen in the following result.

\begin{corollary}
\label{k_unif_all_ints}
Let $H$ be a connected, $k$-uniform hypergraph with $n$ vertices and multiple distinct edges. Apart from $P_n$ and $Q_n$, the partitions in the (lazy) burning distribution are all of the form $\left(\frac{m}{k},\frac{m+1}{k} \right]$ for some $m\in\{0,1,\ldots,k-2\}$, and $Q_n=P_n=\left(\frac{k-1}{k},1\right)$.
\end{corollary}

\begin{proof}
First, $Q_n=P_n=\left(\frac{k-1}{k},1\right)$ due to Lemmas \ref{k_unif_Q1_Qn} and \ref{k_unif_Pn}. 

Now, for $m\in\{0,1,\ldots,k-2\}$, observe that if $p_1,p_2\in \left(\frac{m}{k},\frac{m+1}{k} \right]$ then $b_{L,p_1}(H)=b_{L,p_2}(H)$ and $b_{p_1}(H)=b_{p_2}(H)$, since $m+1$ vertices on fire in an edge cause the fire to propagate whether we use proportion $p_1$ or $p_2$. Hence, for each $m\in\{0,1,\ldots,k-2\}$, the interval $\left(\frac{m}{k},\frac{m+1}{k} \right]$ is a subset of some interval in the (lazy) burning distribution. 

Choose any $m\in\{0,1,\ldots,k-2\}$, and let $\left(\frac{m}{k},\frac{m+1}{k} \right]\subseteq Q_r$ and $\left(\frac{m}{k},\frac{m+1}{k} \right]\subseteq P_s$ for some $r,s<n$. So burning using any proportion $p\in \left(\frac{m}{k},\frac{m+1}{k} \right]$ yields lazy burning number $r$ and burning number $s$. Suppose $Q_r\nsubseteq \left(\frac{m}{k},\frac{m+1}{k} \right]$. Then there is some proportion $p^\prime\in Q_r$ that is not in $\left(\frac{m}{k},\frac{m+1}{k} \right]$. But by Corollary \ref{k_unif_cutoff_cor}, if $p^\prime\leq\frac{m}{k}$ then $b_{L,p^\prime}(H)<r$, so $p^\prime\notin Q_r$. Similarly, if $p^\prime>\frac{m+1}{k}$ then $r<b_{L,p^\prime}(H)$, so $p^\prime\notin Q_r$. These are both contradictions, so our latest supposition must be false. Therefore, $Q_r\subseteq \left(\frac{m}{k},\frac{m+1}{k} \right]$. A similar argument yields $P_s\subseteq \left(\frac{m}{k},\frac{m+1}{k} \right]$. 

We have therefore proven that, given $m\in\{0,1,\ldots,k-2\}$, the interval $\left(\frac{m}{k},\frac{m+1}{k} \right]$ is exactly equal to some interval in the (lazy) burning distribution.
\end{proof}

Given a hypergraph $H$, the set of (lazy) burning numbers of $H$ as $p$ ranges over $(0,1)$ is a subset of $\{1, \ldots, |V(H)|\}$.
Indeed, if $H$ is $k$-uniform, then some values in $\{1, \ldots, |V(H)|\}$ will not be realized for any $p$, as the following result shows.  

\begin{corollary}
\label{n_minus_k_unif}
Let $H$ be a simple, connected, $k$-uniform hypergraph with $n$ vertices and more than one edge. Then there are exactly $n-k$ empty intervals in the (lazy) burning distribution.
\end{corollary}

\begin{proof}
By Corollary \ref{k_unif_all_ints}, the intervals in the (lazy) burning distribution are exactly $\left(0,\frac{1}{k}\right]$, $\left(\frac{1}{k},\frac{2}{k}\right]$, $\ldots$ , $\left(\frac{k-2}{k},\frac{k-1}{k}\right]$, and $\left(\frac{k-1}{k},1\right)$, so there are exactly $k$ nonempty intervals in the (lazy) burning distribution. But there are $n$ possible values for $b_{L,p}(H)$ and $b_p(H)$, so there are $n-k$ empty intervals in the (lazy) burning distribution. 
\end{proof}

Due to Corollaries \ref{k_unif_all_ints} and \ref{n_minus_k_unif}, if $H$ is a simple, connected, $k$-uniform hypergraph with $n$ vertices and more than one edge, then exactly $k$ different values in $\{1,2,\ldots,n\}$ are admissible (lazy) burning numbers for $H$ for some proportion $p$. Indeed, $Q_1$ is nonempty, $P_1$ is empty, and both $Q_n$ and $P_n$ are nonempty by Lemmas \ref{k_unif_Q1_Qn} and \ref{k_unif_Pn}.  

We now consider the burning distribution of possibly non-uniform hypergraphs.

\begin{theorem}
\label{exact_intervals}
Let $H$ be a hypergraph and $\{x_1,x_2,\ldots,x_\ell\}$ be the set of different edge sizes in $H$. Define $X=\left\{ \frac{1}{x_i},\frac{2}{x_i},\ldots,\frac{x_i-1}{x_i}\mid 1\leq i\leq\ell\right\}$, and write its elements in increasing order as $y_1,y_2,\ldots,y_m$. Then each of $(0,y_1],$ $(y_1,y_2],$ $\ldots$ , $(y_{m-1},y_m],$ and $(y_m,1)$ is a subset of some interval in the (lazy) burning distribution.
\end{theorem}

\begin{proof}
By the construction of $X$, if we choose $p$ and $q$ both in one of the intervals $(0,y_1]$, $(y_1,y_2]$, $\ldots$ , $(y_{m-1},y_m]$, or $(y_m,1)$, then $p$ and $q$ are ``close'' enough so that $\lceil p|e|\rceil=\lceil q|e|\rceil$ for every edge $e\in E(H)$. Therefore, $b_{L,p}(H)=b_{L,q}(H)$ and $b_p(H)=b_q(H)$, so the result follows.
\end{proof}

\begin{corollary}
\label{max_num_of_intervals}
Let $H$ be a hypergraph and $\{x_1,x_2,\ldots,x_\ell\}$ be the set of different edge sizes in $H$. Then the maximum number of intervals in the (lazy) burning distribution is $$1-\ell+\sum_{i=1}^\ell x_i.$$
\end{corollary}
 
\begin{proof}
Construct $X$ in the same way as Theorem \ref{exact_intervals}, and write its elements in increasing order as $y_1,y_2,\ldots,y_m$, so $m=|X|$. Observe that $X$ is as large as possible exactly when all of the $x_i$ are relatively prime. In this case, whenever $i\neq j$ we have that $\left\{\frac{1}{x_i},\frac{2}{x_i},\ldots,\frac{x_i-1}{x_i}\right\}$ and $\left\{\frac{1}{x_j},\frac{2}{x_j},\ldots,\frac{x_j-1}{x_j}\right\}$ share no elements. Hence, when $X$ is as large as possible, we have $$m=\sum_{i=1}^{\ell}\left(x_i-1\right)=\left(\sum_{i=1}^{\ell}x_i\right)-\ell.$$ Now, observe that $(0,y_1]$, $(y_1,y_2]$, $\ldots$ , $(y_{m-1},y_m]$, and $(y_m,1)$ are $m+1$ intervals that partition $(0,1)$, and by Theorem \ref{exact_intervals}, each of these intervals is a subset of an interval in the (lazy) burning distribution. Hence, the (lazy) burning distribution contains the maximum number of intervals when $X$ is as large as possible, and each of the intervals $(0,y_1]$, $(y_1,y_2]$, $\ldots$ , $(y_{m-1},y_m]$, and $(y_m,1)$ induces a different (lazy) burning number. In this case, the number of intervals in the (lazy) burning distribution is $$m+1=1-\ell+\sum_{i=1}^\ell x_i.\qedhere$$
\end{proof}

Theorem \ref{exact_intervals} is useful because, given a hypergraph $H$, it tells us which (finitely many) proportions we need to consider in order to deduce the (lazy) burning distribution of $H$. In particular, we only need to compute $b_{y_i}(H)$ for each $i\in\{1,2,\ldots,m\}$, and $b_{p}(H)$ for some $p\in(y_m,1)$ (and similarly for the lazy case). If, for example, $b_{y_i}(H)=b_{y_{i+1}}(H)$, then we know $(y_{i-1},y_i]\cup(y_i,y_{i+1}]$ is a subset of an interval in the burning distribution. Otherwise, if $b_{y_i}(H)<b_{y_{i+1}}(H)$, then $y_i$ is the right endpoint of an interval in the burning distribution.

In light of Corollary \ref{max_num_of_intervals}, it seems as though a hypergraph with many different edge sizes that are relatively prime is likely to have more intervals in its (lazy) burning distribution. We leave this as an open problem.

We previously conjectured the following: if $p,q\in(0,1)$ and $p\leq q$, then for any hypergraph $H$, $b_q(H)-b_{L,q}(H)\leq b_p(H)-b_{L,p}(H)$. Informally, this means that as the proportion increases, the difference between the burning number and lazy burning number decreases monotonically. However, this was disproven after obtaining our computational results. Observe Table \ref{burning_dist_table_1} on page \pageref{burning_dist_table_1}, and notice that there are several counterexamples, one being the BIBD$(13,4,1)$. The following weaker statement is still open.

\begin{conjecture}
\label{limit_conj_special_case} 
Let $H$ be a hypergraph, $p,q\in(0,1)$, and $b_{L,p}(H)=b_{L,q}(H)$. Then $b_p(H)=b_q(H)$.
\end{conjecture}

\begin{lemma}
\label{lazy_subset_round_based}
Let $H$ be a hypergraph. Then, every interval in the lazy burning distribution is a subset of some interval in the burning distribution if and only if $b_{L,p}(H)=b_{L,q}(H)$ implies $b_p(H)=b_q(H)$ for any $p,q\in(0,1)$.
\end{lemma}

\begin{proof}
For the forward implication, assume every interval in the lazy burning distribution is a subset of some interval in the burning distribution, and choose $p,q\in(0,1)$ such that $b_{L,p}(H)=b_{L,q}(H)$. Then $p$ and $q$ belong to the same interval in the lazy burning distribution. Since this interval is a subset of an interval in the burning distribution, $p$ and $q$ also belong to the same interval in the burning distribution. Hence, $b_p(H)=b_q(H)$.

Now, we prove the contrapositive of the reverse implication. Assume there exists an interval $I_0$ in the lazy burning distribution that is not a subset of any interval in the burning distribution. Then, there must exist two intervals $I_1$ and $I_2$ in the burning distribution such that $I_0\cap I_1$ and $I_0\cap I_2$ are both nonempty. Less formally, $I_0$ ``overlaps'' with multiple intervals in the burning distribution. Choose $p$ and $q$ such that $p\in I_0\cap I_1$ and $q\in I_0\cap I_2$. Then $p$ and $q$ are in the same interval in the lazy burning distribution, but different intervals in the burning distribution. Therefore, we have found $p,q\in(0,1)$ such that $b_{L,p}(H)=b_{L,q}(H)$ and $b_p(H)\neq b_q(H)$.
\end{proof}

Due to Lemma \ref{lazy_subset_round_based}, Conjecture \ref{limit_conj_special_case} may be re-formulated as follows: \emph{for any hypergraph $H$, every interval in the lazy burning distribution is a subset of some interval in the burning distribution}. Indeed, this is true for all $k$-uniform hypergraphs $H$; see Corollary \ref{k_unif_intsss}.

\begin{corollary}
\label{k_unif_intsss}
Let $H$ be a $k$-uniform hypergraph and $p,q\in(0,1)$. Then, $b_{L,p}(H)=b_{L,q}(H)$ if and only if $b_p(H)=b_q(H)$.
\end{corollary}

\begin{proof}
By Corollary \ref{k_unif_all_ints}, the intervals in the lazy burning distribution are identical to those in the burning distribution (they both contain exactly $\left(0,\frac{1}{k}\right], \left(\frac{1}{k},\frac{2}{k}\right],\ldots$ , and $\left(\frac{k-1}{k},1\right)$).
\end{proof}

To disprove Conjecture  \ref{limit_conj_special_case}, it would suffice to find a hypergraph $H$ in which some interval in the lazy burning distribution is not a subset of any interval in the burning distribution.

\section{Computational Results}
\label{comp_section}

In this section we obtain computational results on the lazy burning distribution of a special family of hypergraphs known as \emph{balanced incomplete block designs}. We chose to investigate this family due to their rigorously studied structure, and indeed we were able to uncover a relationship between the lazy burning process and the design-theoretic properties of a balanced incomplete block design. 

We use \cite{design_text} as a reference for the theory of designs. A \emph{design} is a pair $(X, \cal{A})$ where $X$ is a set of elements called \emph{points} and $\cal{A}$ is a multiset consisting of nonempty subsets of $X$, called \emph{blocks}. Now, let $v,k,\lambda \in\mathbb{N}$ such that $v>k\geq 2$. A $(v,k,\lambda)$-\emph{balanced incomplete block design}, abbreviated \emph{BIBD}$(v,k,\lambda)$, is a design $(X,\cal{A})$ such that  $|X|=v$, each block contains exactly $k$ points, and every pair of distinct points is contained in exactly $\lambda$ blocks. Observe that a design $(X,\cal{A})$ induces a hypergraph $H$ on which the burning game can be played, since we can take $V(H)=X$ and $E(H)=\cal{A}$. 

We wrote a C program that computes the (lazy) burning distribution of a hypergraph, and tested it on many of the small balanced incomplete block designs found in \cite{designs_handbook}; see Table \ref{burning_dist_table_1}. We use the same numbering as in \cite{designs_handbook} when there are multiple non-isomorphic designs with the same parameters. For the sake of readability, we write the (lazy) burning distributions in a condensed format. The list of numbers in the (lazy) burning distribution column represent the (lazy) burning numbers of the hypergraph that are realized as $p$ ranges over $(0,1)$. Then, one may apply Theorem \ref{comparing_proportions} and Corollary \ref{k_unif_all_ints} in order to determine the (lazy) burning distribution. For example, the intervals in the (lazy) burning distribution for the BIBD$(16,4,1)$ are exactly $\left(0,\frac{1}{4}\right]$, $\left(\frac{1}{4},\frac{1}{2}\right]$, $\left(\frac{1}{2},\frac{3}{4}\right]$, and $\left(\frac{3}{4},1\right]$. Then, using the information in the table, we must have $Q_1=\left(0,\frac{1}{4}\right]$, $Q_3=\left(\frac{1}{4},\frac{1}{2}\right]$, $Q_7=\left(\frac{1}{2},\frac{3}{4}\right]$, $Q_{16}=\left(\frac{3}{4},1\right]$, $P_1=\left(0,\frac{1}{4}\right]$, $P_5=\left(\frac{1}{4},\frac{1}{2}\right]$, $P_8=\left(\frac{1}{2},\frac{3}{4}\right]$, and $P_{16}=\left(\frac{3}{4},1\right]$. 

When testing the BIBD$(25,4,1)$s in the round-based game with proportion $p=\frac{3}{4}$, the program ran for too long and would not finish in a reasonable amount of time. However, we were able to find a burning sequence of length 9 for each BIBD$(25,4,1)$, so their burning numbers are at most 9 for any proportion in $\left(\frac{1}{2},\frac{3}{4}\right]$. Observe that the last two rows of Table~\ref{burning_dist_table_1} have a ``$\leq 9$'' in the corresponding position to reflect this fact. In light of Theorem \ref{tight_ineq}, this value could be either 8 or 9 for the BIBD$(25,4,1)$s \#1 and \#6, and for the rest it could be 7, 8, or 9.

\begin{table}[]
\centering
\begin{tabular}{|c|c|c|}
\hline
BIBD$(v,k,\lambda)$ & Lazy burning distribution & Burning distribution \\ \hline
$(6,3,2)$           & 1, 2, 6                   & 2, 4, 6              \\ \hline
$(7,3,1)$           & 1, 3, 7                   & 2, 4, 7              \\ \hline
$(7,3,2)$ \#1       & 1, 3, 7                   & 2, 4, 7              \\ \hline
$(7,3,2)$ \#2 - 4     & 1, 2, 7                   & 2, 4, 7              \\ \hline
$(7,3,3)$ \#1       & 1, 3, 7                   & 2, 4, 7              \\ \hline
$(7,3,3)$ \#2 - 10    & 1, 2, 7                   & 2, 4, 7              \\ \hline
$(8,4,3)$ \#1 - 3     & 1, 2, 3, 8                & 2, 3, 5, 8           \\ \hline
$(8,4,3)$ \#4       & 1, 2, 4, 8                & 2, 3, 5, 8           \\ \hline
$(9,3,1)$           & 1, 3, 9                   & 2, 5, 9              \\ \hline
$(9,3,2)$ \#1       & 1, 3, 9                   & 2, 5, 9              \\ \hline
$(9,3,2)$ \#2 - 36    & 1, 2, 9                   & 2, 4, 9              \\ \hline
$(9,4,3)$ \#1 - 10    & 1, 2, 3, 9                & 2, 3, 5, 9           \\ \hline
$(9,4,3)$ \#11      & 1, 2, 4, 9                & 2, 3, 5, 9           \\ \hline
$(10,4,2)$ \#1 - 2    & 1, 2, 4, 10               & 2, 4, 6, 10          \\ \hline
$(10,4,2)$ \#3      & 1, 2, 5, 10               & 2, 4, 6, 10          \\ \hline
$(11,5,2)$          & 1, 2, 4, 6, 11            & 2, 4, 5, 8, 11       \\ \hline
$(13,4,1)$          & 1, 3, 6, 13               & 2, 4, 8, 13          \\ \hline
$(13,3,1)$ \#1 \& 2    & 1, 3, 13                  & 2, 5, 13             \\ \hline
$(15,3,1)$ \#1      & 1, 4, 15                  & 2, 5, 15             \\ \hline
$(15,3,1)$ \#2 - 80   & 1, 3, 15                  & 2, 5, 15             \\ \hline
$(15,7,3)$ \#1           & 1, 2, 3, 5, 7, 10, 15           & 2, 3, 4, 6, 8, 11, 15          \\ \hline
$(15,7,3)$ \#2 - 5         & 1, 2, 3, 5, 7, 9, 15            & 2, 3, 4, 6, 8, 11, 15        \\ \hline 
$(16,4,1)$          & 1, 3, 7, 16               & 2, 5, 8, 16          \\ \hline
$(16,6,2)$ \#1           & 1, 2, 4, 6, 10, 16             & 2, 4, 5, 8, 11, 16            \\ \hline
$(16,6,2)$ \#2 \& 3         & 1, 2, 4, 6, 9, 16              & 2, 4, 5, 8, 11, 16           \\ \hline
$(19,9,4)$ \#1 - 6         & 1, 2, 3, 4, 6, 8, 10, 13, 19      & 2, 3, 4, 6, 7, 9, 11, 14, 19                    \\ \hline
$(25,4,1)$ \#1 \& 6      & 1, 3, 7, 25                  &     2, 5, $\leq 9$, 25                 \\ \hline
$(25,4,1)$ \#2 - 5 \& 7 - 18 & 1, 3, 6, 25                  &    2, 5, $\leq 9$, 25                   \\ \hline
\end{tabular}
\caption{Burning and lazy burning distributions for some small BIBDs}
\label{burning_dist_table_1}
\end{table}

Recall that playing the proportion-based lazy burning game on a $k$-uniform hypergraph with $p=\frac{k-1}{k}$ is equivalent to playing the \emph{original} lazy burning game on that hypergraph, since then the two propagation rules coincide. There are many instances in Table \ref{burning_dist_table_1} in which one of the BIBDs has a higher lazy burning number than all the other non-isomorphic BIBDs with the same parameters, and indeed this only occurs when $p=\frac{k-1}{k}$. Since we are observing non-isomorphic BIBDs with the same parameters that have different lazy burning numbers, it is natural to ask what other properties these BIBDs have, and if any of them correlate with having a higher lazy burning number. 

Before moving on, let us recall some definitions. Given two hypergraphs $H$ and $G$, an \emph{isomorphism} is a bijective function $f: V(H)\rightarrow V(G)$ with the property that $e$ is an edge in $H$ if and only if $\{f(v)\mid v\in e\}$ is an edge in $G$. An isomorphism $f:V(H)\rightarrow V(H)$ is called an \emph{automorphism} of $H$. A hypergraph $H$ is \emph{vertex-transitive} if for every pair of vertices $u,v\in V(H)$, there exists an automorphism $f$ of $H$ with $f(u)=v$. It is well-known that the set of all automorphisms of a finite hypergraph forms a group under composition of functions. The \emph{automorphism group order} of a hypergraph $H$ is the number of distinct automorphisms that exist for $H$. 

We implemented Nauty and Traces \cite{Nauty_Traces} in a C program and computed the automorphism group orders of the relevant hypergraphs in Table \ref{burning_dist_table_1} -- these results can be seen in Table \ref{auto_group_table}. Indeed every time a BIBD had a higher lazy burning number than the other non-isomorphic BIBDs with the same parameters, it also had the highest automorphism group order among all the BIBDs with those parameters. It is therefore reasonable to suspect that there is a strong correlation between the automorphism group order and the lazy burning number (in the original game) of a BIBD. The only outliers are the BIBD$(25,4,1)$s, as there are two (\#1 and \#6) with lazy burning number seven, while the rest have lazy burning number six. However, the BIBD$(25,4,1)$s still support the hypothesis that there is a connection between automorphism group order and the lazy burning number, since numbers 1 and 6 have the first and second highest automorphism group orders respectively. 

Initially, our intuition on the connection between the automorphism group order and (lazy) burning number of a hypergraph $H$ was as follows. A high automorphism group order means that $H$ has many automorphisms, and hence $H$ is more likely to be vertex-transitive. It was recently proven in \cite{edge_conn_vert_trans} that a vertex-transitive hypergraph is maximally edge-connected (\emph{if} it is also connected, linear, and uniform). A maximally edge-connected hypergraph has a high ``density'' of edges in some sense, which should expedite the spread of the fire and result in lower values for $b_L(H)$ and $b(H)$. Hence, we expected that a higher automorphism group order would correlate with a lower (lazy) burning number (at least for the BIBDs with $\lambda=1$, as these are all linear). However, the results from Table \ref{auto_group_table} do not support this intuition, and in fact suggest that the opposite correlation may be true. It would be interesting to determine the flaw with this intuition, as this may lead us to a more fundamental understanding of (lazy) hypergraph burning.


\begin{table}[]
\centering
\begin{tabular}{|c|c|c|}
\hline
BIBD$(v,k,\lambda)$ &  Lazy Burning Number &Automorphism Group Order\\             \hline
$(7,3,2)$ \#1             &  3  &        $21504$                 \\  \hline
$(7,3,2)$ \#2 - 4         &  2 &        $\leq 384$               \\ \hline
$(7,3,3)$ \#1              &    3   &       $47029248$             \\ \hline
$(7,3,3)$ \#2 - 10        &  2    &       $\leq 82944$            \\ \hline
$(8,4,3)$ \#1 - 3          &  3 &       $\leq 48$                  \\ \hline
$(8,4,3)$ \#4               &  4 &       $1344$                      \\ \hline
$(9,3,2)$ \#1                &   3   &       $1769472$               \\ \hline
$(9,3,2)$ \#2 - 36          &  2   &       $\leq 1536$              \\ \hline
$(9,4,3)$ \#1 - 10         & 3  &       $\leq 32$                  \\ \hline
$(9,4,3)$ \#11               &  4 &      $144$                       \\ \hline
$(10,4,2)$ \#1 - 2           &  4   &      $\leq 48$                \\ \hline
$(10,4,2)$ \#3                &   5  &      $720$                 \\ \hline
$(15,3,1)$ \#1               &  4 &     $20160$           \\ \hline
$(15,3,1)$ \#2 - 80        &  3 &     $\leq 288$         \\ \hline
$(15,7,3)$ \#1               &   10    &     $40320$            \\ \hline
$(15,7,3)$ \#2 - 5           &   9   &     $\leq 1152$        \\ \hline 
$(16,6,2)$ \#1                &  10 &     $23040$            \\ \hline
$(16,6,2)$ \#2 \& 3         & 9  &     $\leq 1536$       \\ \hline
$(25,4,1)$ \#1               & 7 &    $504$                  \\ \hline
$(25,4,1)$ \# 6            &   7 &      $150$                \\ \hline
$(25,4,1)$ \#2 - 5 \& 7 - 18  & 6 &    $\leq 63$  \\ \hline
\end{tabular}
\caption{Lazy burning numbers automorphism group orders  for some small BIBDs}
\label{auto_group_table}
\end{table}

\section{Summary and Open Problems}

In Section \ref{general_results} we proved that if the proportion is of the form $p=\frac{1}{n}$ for some $n\in\mathbb{N}\setminus\{1\}$ and $H$ is connected, then the lazy burning number of $H$ with respect to $p$ is bounded above by $\left\lceil \frac{|V(H)|}{n}\right\rceil$; see Theorem \ref{pproptheorem}. It is natural to ask if the analogous bound holds when $p$ is \emph{not} of the form $\frac{1}{n}$; see Conjecture \ref{ppropconj2}. In fact, we have shown that such a bound holds for all $p\in(0,1)$ if and only if it holds for all rational $p$ in $(0,1)$ not of the form $\frac{1}{n}$, so one may assume $p$ is rational when attempting a proof. Note that the techniques used in the proof of Theorem \ref{pproptheorem} may not work when attempting to prove Conjecture \ref{ppropconj2}, as is discussed shortly after the conjecture. Also, the conjecture stipulates that each edge in $H$ is flammable, which is always true when $p$ is of the form $\frac{1}{n}$. Do we need to assume this in the statement of the conjecture?

In Corollary \ref{max_num_of_intervals} we calculate the maximum possible number of intervals in the (lazy) distribution of a hypergraph $H$. The proof of this theorem  shows that this maximum value can be attained only if the sizes of all the edges in $H$ are relatively prime. Can more be said about how the number of relatively prime edge sizes in $H$ influences the number of intervals in the (lazy) burning distribution? It seems very likely that $k$-uniform hypergraphs have the least possible number of intervals in their (lazy) burning distribution. Can this be proven?

We were able to determine exactly what intervals are in the (lazy) burning distribution of a $k$-uniform hypergraph, but a method to determine the corresponding (lazy) burning numbers still eludes us. Can we determine such a method for BIBDs, or other $k$-uniform designs with a rigorously defined structure?

Through computational methods, we were able to disprove the following statement: if $p\leq q$ then $b_q(H)-b_{L,q}(H)\leq b_p(H)-b_{L,p}(H)$. Informally, this statement would mean that as the proportion $p$ increases, the difference between $b_p(H)$ and $b_{L,p}(H)$ decreases monotonically. The following weaker statement is still open: if $b_{L,p}(H)=b_{L,q}(H)$ then $b_p(H)=b_q(H)$ -- see Conjecture \ref{limit_conj_special_case}. We were able to show that this statement is true if and only if every interval in the lazy burning distribution is a subset of an interval in the burning distribution. Hence, if one wishes to disprove the conjecture, it would suffice to find a hypergraph $H$ with an interval in its lazy burning distribution that is not contained in any interval in its burning distribution.

The computational results from Section \ref{comp_section} make a very strong case for a relationship between the automorphism group order of a BIBD and its lazy burning number (from the original model). It may be true that simply knowing the automorphism group order of a BIBD $H$ is not enough to infer anything about $b_L(H)$. However, in light of the results from Table \ref{auto_group_table}, the following conjecture seems plausible.

\begin{conjecture}
If $H$ and $G$ are two BIBDs on the same parameters and $H$ has a higher automorphism group order than $G$, then $b_L(H)\geq b_L(G)$. 
\end{conjecture}

Finally, many of the open problems presented at the end of \cite{our_paper!} can equally apply to our proportion-based variant of hypergraph burning.

\section{Acknowledgements}

Authors Burgess and Pike acknowledge NSERC Discovery Grant support, and Jones acknowledges NSERC and AARMS scholarship support. Jones also acknowledges the support of MITACS, which partially funded this research through a collaborative grant with Verafin. We also thank Verafin for both their financial support
and the contributions of their researchers who co-authored this paper.















\begin{thebibliography}{99}

\addcontentsline{toc}{section}{References}

\bibitem{VERY_first_ref} N. Alon. Transmitting in the $n$-Dimensional Cube. \emph{Discrete Applied Mathematics} Vol. 37-38 (1992) pp. 9-11. 

\bibitem{bahm_sanj} M. Bahmanian, M. \v{S}ajna. Connection and Separation in Hypergraphs. \emph{Theory and Applications of Graphs} Vol. 2, Iss. 2, Art. 5 (2015).

\bibitem{bootstrap_paper_1} J. Balogh, B. Bollob\'{a}s, R. Morris, O. Riordan. Linear Algebra and Bootstrap Percolation. \emph{Journal of Combinatorial Theory Series A} Vol. 119 (2012) pp. 1328--1335.

\bibitem{new_paper} P. Bastide, M. Bonamy, A. Bonato, P. Charbit, S. Kamali, T. Pierron, M. Rabie. Improved Pyrotechnics: Closer to the Burning Number Conjecture. arXiv:2110.10530v2 (2022)

\bibitem{burning_is_hard} S. Bessy, A. Bonato, J. Janssen, D. Rautenbach, E. Roshanbin. Burning a Graph is Hard. \emph{Discrete Applied Mathematics} Vol. 232 (2017) pp. 73-87. 

\bibitem{bounds} S. Bessy, A. Bonato, J. Janssen, D. Rautenbach, E. Roshanbin. Bounds on the Burning Number. \emph{Discrete Applied Mathematics} Vol. 235 (2018) pp. 16-22. 

\bibitem{Bonato_summary} A. Bonato. A Survey of Graph Burning. \emph{Contributions to Discrete Mathematics} Vol. 16 No. 1 (2021) pp. 185-197.

\bibitem{first_paper} A. Bonato, J. Janssen, E. Roshanbin. Burning a Graph as a Model of Social Contagion. Proceedings of WAW (2014) pp. 13-22.

\bibitem{how_to} A. Bonato, J. Janssen, E. Roshanbin. How to Burn a Graph. \emph{Internet Mathematics} 1-2 (2016) pp. 85-100.

\bibitem{our_Latin_paper} A. Bonato, C. Jones, T. Marbach, T. Mishura. How to Burn a Latin Square. (Preprint, 2024)

\bibitem{hypergraph_text} A. Bretto. Hypergraph Theory an Introduction. \emph{Springer} (2013).

\bibitem{zero_forcing_complexity} B. Brimkov, C. Fast, I. Hicks. Computational Approaches for Zero Forcing and Related Problems. \emph{European Journal of Operational Research} 273 (2019) pp. 889-903.

\bibitem{our_STS_paper} A. Burgess, P. Danziger, C. Jones, T. Marbach, D. Pike. Burning Steiner Triple Systems. (Preprint, 2024).

\bibitem{our_paper!} A. Burgess, C. Jones, D. Pike. Extending Graph Burning to Hypergraphs. (Preprint, 2024).

\bibitem{edge_conn_vert_trans} A. Burgess, R. Luther, D. Pike. The Edge-Connectivity of Vertex-Transitive Hypergraphs. \emph{Journal of Graph Theory} Vol. 105 (2) (2023) pp. 252-259.

\bibitem{designs_handbook} C. Colbourn and J. Dinitz. Handbook of Combinatorial Designs 2nd ed. \emph{Chapman \& Hall/CRC} (2007).

\bibitem{triple_systems_text} C. Colbourn and A. Rosa. Triple Systems. \emph{Clarendon, Oxford University Press} (1999).

\bibitem{conn_in_hyp} M. Dewar, D. Pike, J. Proos. Connectivity in Hypergraphs. \emph{Canad. Math. Bull.} Vol. 61 (2) (2018) pp. 252-271.

\bibitem{STS_with_no_subsystems} J. Doyen. Sur la structure de certaines syst\`{e}mes triples de Steiner. \emph{Math. Z.} 111 (1969) pp. 289-300.

\bibitem{mythesis} C. Jones. Hypergraph Burning. Masters Thesis, Memorial University of Newfoundland, 2023.

\bibitem{all_sts_15s} R. Mathon, K. Phelps, A. Rosa. Small Steiner Triple Systems and their Properties. \emph{Ars Combinatoria} Vol. 15 (1983) pp. 3-110 (and errata: Vol 16 (1983) p. 286).

\bibitem{Nauty_Traces} B. McKay and A. Piperno. Practical Graph Isomorphism, \textrm{II}. \emph{Journal of Symbolic Computation} Vol. 60 (2014) pp. 94-112.

\bibitem{asymptotically_paper} S. Norin and J. Turcotte. The Burning Number Conjecture Holds Asymptotically. 
{\em Journal of Combinatorial Theorey Series B} Vol. 168 (2024), 208--235.


\bibitem{Singer} J. Singer. A Theorem in Finite Projective Geometry and Some Applications to Number Theory. \emph{Transactions of the American Mathematical Society} Vol. 43 No. 3 (1938) pp. 377-385.

\bibitem{design_text} D. Stinson. Combinatorial Designs: Constructions and Analysis. \emph{Springer} (2004).

\bibitem{graph_theory_west} D. West. Introduction to Graph Theory 2nd ed. \emph{Pearson} (2001).

\end{thebibliography}
\end{document}